\documentclass[12pt]{article}
\usepackage[expansion=false]{microtype}
\usepackage[reqno]{amsmath}
\usepackage{fixmath}
\usepackage{amssymb, amsthm, enumitem}
\usepackage{mathtools}
\usepackage{graphicx}

\usepackage{ifthen,amssymb,amsfonts,amsmath,amsthm,graphicx, color, tikz, pgf}

\usepackage{hyperref}

\newtheoremstyle{plainsl}%
    {\topsep}
    {\topsep}
    {\slshape} 
    {}
    {\normalfont\bfseries}
    {.}
    { }
    {}

\swapnumbers

{\theoremstyle{plainsl}
\newtheorem{thm}{Theorem}[section]
\newtheorem{lem}[thm]{Lemma}
\newtheorem{cor}[thm]{Corollary}
}
{\theoremstyle{remark}
  }

\renewcommand\proof{\noindent\textsl{Proof. }}
\newcommand\sqr[2]{{\vbox{\hrule height.#2pt
   \hbox{\vrule width.#2pt height#1pt \kern#1pt
        \vrule width.#2pt}\hrule height.#2pt}}}
\renewcommand\qed{%
    \ifmmode\eqno\sqr53
    \else\nolinebreak\ \hfill\sqr53\medbreak\fi}

\numberwithin{equation}{section}

\newcommand{\abs}[1]{\left | #1 \right |}
\newcommand{\lp}{\left (}
\newcommand{\rp}{\right )} \newcommand{\lb}{\left [}
  \newcommand{\rb}{\right ]} \newcommand{\lsb}{\left \{ }
\newcommand{\rsb}{\right \} }

\newcommand\epvxs[2]{\epsilon_{#1,#2}}
\newcommand\eab{\epvxs{a}{b}}
\newcommand\comp[1]{{\mkern2mu\overline{\mkern-2mu#1}}}

\newcommand{\ver}[1]{\mathbf{e}_{ #1}}

\newcommand\cx{\mathbb{C}}

\title{Perfect State Transfer on Oriented Graphs}

\author{Chris Godsil and Sabrina Lato}

\begin{document}
\maketitle

\begin{abstract}
  Quantum walks on undirected graphs have been studied using symmetric matrices, such as the adjacency or Laplacian matrix, and many results about perfect state transfer are known. We extend some of those results to oriented graphs. We also study the phenomena, unique to oriented graphs, of multiple state transfer, where there is a set of vertices such that perfect state transfer occurs between every pair in that set. We give a characterization of multiple state transfer, and a new example of a graph where it occurs.\end{abstract}

\section{Introduction}

A continuous quantum walk on a graph is defined by taking some Hermitian matrix \( H \) and considering the time-dependent unitary matrix
\[ U (t) = \exp(-i t H). \]
We call \(U\) the \textit{transition matrix} of the walk. In most previous work, $H$ is chosen to be the 
adjacency matrix or the Laplacian matrix of an undirected graph.

Cameron et al \cite{Cameron2014a} considered \textit{oriented graphs} where every edge is given a unique direction.
We take the adjacency matrix of an oriented graph to be the matrix \(A\) with rows and columns indexed
by the vertices of the graph, and \(A_{a,b}\) equal to \(1\) if the the edge \(\{a,b\}\) is
oriented from \(a\) to \(b\), equal to \(-1\) if the the edge \(\{a,b\}\) is
oriented from \(a\) to \(b\), and equal to zero if \(a\) and \(b\) are not adjacent. Consequently
\(A\) is skew symmetric. If \(A\) is a skew symmetric matrix, then \(iA\) is Hermitian and so
\[
	U(t) = \exp(-it(iA)) = \exp(tA)
\]
is the transition matrix of a continuous quantum walk. Although \(U(t)\) is real and orthogonal, this defines a continuous quantum walk, not a classical random walk.

Let \( \ver{a} \) denote the characteristic vector for the vertex \( a \). We say that \textit{perfect state transfer} occurs from vertex \( a \) to vertex \( b \) if there is some time \( \tau \in \mathbb{R} \) and some complex phase factor \( \lambda \) with norm one such that
\[ U (\tau) \ver{a} = \lambda \ver{b}. \]
This definition matches that used for unoriented graphs. Our goal in this paper is to characterize
the case where perfect state transfer occurs on an oriented graph. For undirected graphs (with the
adjacency matrix or the Laplacian as Hamiltonian), Kay \cite{Kay} showed that for any vertex
\(a\) in a graph, there is most one vertex \(b\) such that we have perfect state transfer
from \(a\) to \(b\). Cameron et al gave two examples of oriented graphs with what they called
\textit{universal state transfer}, oriented graphs such that for each pair of vertices \(a\)
and \(b\), there was perfect state transfer from \(a\) to \(b\) at some time. Here we provide
an oriented graph on eight vertices with a subset \(C\) of four vertices, such that there is
perfect state transfer between any pair of vertices in \(C\). We call this 
\textit{multiple state transfer}.

\section{Basic Properties of Quantum Walks}

Skew symmetric matrices have some convenient properties, which are standard results from linear algebra and can be found in texts such as Zhang~\cite{Zhang1999}. Among them, the eigenvalues have no real part and are symmetric about the real axis. Skew symmetric matrices are also normal, and therefore admit a spectral decomposition as described
(for example) in Godsil and Royle~\cite{Godsil2001}. That is, for any skew symmetric matrix \( A \) with distinct eigenvalues \( \theta_0, \ldots, \theta_d \) we may write
\[ 
	A = \sum_{r=0}^d \theta_r E_r 
\]
where the \textit{spectral idempotents} \( E_r \) represent projection into the \( \theta_r \)-eigenspace. 
These spectral idempotents are Hermitian, positive semidefinite, pairwise orthogonal, and sum to the identity. 
For the skew symmetric case, we have the further useful property that if \( \theta_r \) is an eigenvalue 
with spectral idempotent \( E_r, \) then \( - \theta_r \) is an eigenvalue with spectral idempotent \( \overline{E_r} \).

Using spectral decomposition, any function \( f \) that converges on the spectrum of \( A \) can be written
\[ 
	f(A) = \sum_{r=0}^d f(\theta_r) E_r. 
\]
In particular,
\[ 
	U(t) = \sum_{r=0}^d e^{t \theta_r} E_r. 
\]
Details and proofs of these facts are given for symmetric matrices in~\cite{Godsil2001} and for skew symmetric matrices in~\cite{Lato2019}. 

We note that \( U(t)\) is real and that 
\[
	\comp{U(t)} = U(t)^{-1} = U(-t).
\]

For unoriented graphs, little is known about the phase factor. All examples of perfect state transfer have 
a root of unity as the phase factor, but it is not known whether this is always true~\cite{Coutinho2020}. 
For oriented graphs, we know something much stronger.

\begin{lem}\label{phase}
	Let \( X \) be an oriented graph with perfect state transfer. Then the phase factor is \( \pm 1. \)
\end{lem}

\proof 
Suppose there is perfect state transfer from \( a \) to \( b \). Then we have some time \( \tau \) 
and some phase factor \( \lambda \) such that $U(t)\ver{a}=\lambda e_b$; since $U(t)$, $\ver{a}$ and $\ver{b}$ are
real it follows that $\lambda$ is real. As $\abs{\lambda}=1$, the lemma follows.\qed

In subsequent computations with perfect state transfer, we will use \( \eab \) to denote the phase factor,
i.e., if there is perfect state transfer from $a$ to $b$, then we write $U(t)\ver{a}=\eab U(t)\ver{b}$.

\section{Switching Automorphisms}

An automorphism of a graph with adjacency matrix \( A \) can be viewed as a permutation matrix \( P \) such that
\[ A = P^T A P. \]

A \textit{monomial matrix} is the product of a permutation matrix with a diagonal matrix. For our purposes, the diagonal matrix will always have entries \( \pm 1 \). Note that if \( \tilde{P} \) is such a monomial matrix, then \( \tilde{P}^T\tilde{P} = I. \) A graph \( X \) has a \textit{switching automorphism} if there exists a monomial matrix \( \tilde{P} \) such that
\[ A = \tilde{P}^T A \tilde{P}. \]

We have defined switching automorphisms to be matrices (acting on \(\cx^{V(X)}\)), but if \(a\in V(X)\)
and \(\tilde P\) is a switching, 
we will sometimes use \( \tilde{P} a \) to denote the vertex such that \(  \ver{\tilde{P}a} = \tilde{P} \ver{a}. \)

Cameron et al~\cite{Cameron2014a} studied the switching automorphism group of graphs with universal state transfer.
We will go in a different direction, and characterize when perfect state transfer can occur from vertex \( a \) to \( \tilde{P} a \) for some switching automorphism \( \tilde{P}. \) The first point to note is that switching automorphisms preserve the quantum walk.

\begin{lem}\label{switching}
	Let \( X \) be an oriented graph with switching automorphism  \( \tilde{P} \) and transition matrix \( U (t) \). Then
    \[ 
		U(t) = \tilde{P}^T U(t) \tilde{P}. 
	\]
\end{lem}

\proof
Since $\tilde{P}$ commutes with $A$, it commutes with $U=\exp(tA)$.\qed
 
This is particularly relevant when there is a switching automorphism along with perfect state transfer.

\begin{cor}\label{automorphismPST}
	Let \( X \) be an oriented graph with switching automorphism \( \tilde{P} \). If \( X \) has perfect state transfer 
	from vertex \( a \) to vertex \( b \) at time \( \tau, \) then there is perfect state transfer 
	from \( \tilde{P} a \) to \( \tilde{P} b \) at the same time \( \tau. \)
\end{cor}

\proof 
Since there is perfect state transfer,
\[ U(\tau) \ver{a} = \pm \ver{b}. \]
By Lemma~\ref{switching},
\[ \tilde{P}^T U(\tau) \tilde{P} \ver{a} = \pm \ver{b}, \]
so
\[ U(\tau) \tilde{P} \ver{a} = \pm \tilde{P} \ver{b} \]
and there is perfect state transfer from \( \tilde{P} a \) to \( \tilde{P} b \) at time \( \tau.\)\qed

\begin{cor}\label{multipleAuto}
	Let \( X \) be an oriented graph with switching automorphism \( \tilde{P} \) such that there is perfect state transfer from 
	vertex \( a \) to vertex \( \tilde{P} a \) at time \( \tau. \) Then for any positive integer \( n \) 
	there is perfect state transfer from \( a \) to \( \tilde{P}^n a \) at time \( n \tau \).
\end{cor}

\proof 
We will prove this by induction. The base case when \( k = 1 \)  is true by definition, so let \( k \) 
be a positive integer such that there is perfect state transfer from \( a \) to \( \tilde{P}^k a \) 
at time \( k \tau \). Then
\[ U((k+1))  \tau) \ver{a} = U(\tau) U(k \tau) \ver{a}. \]
By the inductive hypothesis there is perfect state transfer from \( a \) to \( \tilde{P}^k a \) 
at time \( k \tau \) and thus
\[ 
	U(\tau) U(k \tau) \ver{a} = \pm U(\tau) \tilde{P}^k \ver{a}. 
\]
Since \( \tilde{P}^k \) is a switching automorphism, by Corollary~\ref{automorphismPST} there is perfect state transfer from \( \tilde{P}^k a \) to \( \tilde{P}^{k+1} a \) at time \( \tau \) and therefore
\[ \pm U(\tau) \tilde{P}^k \ver{a} = \pm \tilde{P}^{k+1} \ver{a}. \]
Thus there is perfect state transfer from \( a \) to \( \tilde{P}^{k+1} a \) at time \((k + 1) \tau \).\qed

\section{Multiple State Transfer}

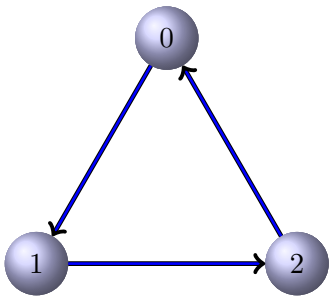
\begin{figure}    
  \begin{center}
    \begin{tikzpicture}[node distance = 2 cm]
      \tikzset{VertexStyle/.style = {shape= circle, ball color= blue!20!white, inner sep= 2pt, outer sep= 0pt, minimum size = 24 pt}}
      \tikzset{EdgeStyle/.style   = {->, double= blue, double distance = 1pt}}
          
      \draw (90:2 cm) node[VertexStyle, font=\small] (0) {0} ;
      \draw (210:2 cm) node[VertexStyle, font=\small] (1) {1} ;
      \draw (330:2 cm) node[VertexStyle, font=\small] (2) {2} ;

      \draw [EdgeStyle] (0)--(1) ;
      \draw [EdgeStyle] (1)--(2) ;
      \draw [EdgeStyle] (2)--(0) ;
      
    \end{tikzpicture}
  \end{center}
  \caption{One orientation of the complete graph on 3 vertices} \label{k3} 
\end{figure}

Kay \cite{Kay} showed that if a real symmetric matrix has perfect state transfer from \( a \) to some other 
vertex \( b, \) it cannot have perfect state transfer from \(a \) to any vertex other than \(a\). This is not the case with oriented graphs, as the complete graph on three vertices oriented as in Figure~\ref{k3} has perfect state transfer from vertex 0 to vertex 1 at time \( \frac{2 \pi}{3 \sqrt{3}} \) and from vertex 0 to 1 at time \( \frac{4 \pi}{3 \sqrt{3}} \).

Cameron et al~\cite{Cameron2014a} studied graphs like this with \textit{universal state transfer}, where there is perfect state transfer between every pair of vertices, although the only known examples are \( K_2 \) and \( K_3 \). Connelly et al~\cite{Connelly} conjectured that these were in fact the only examples of unweighted oriented graphs with universal state transfer.

We will consider the more general case of \textit{multiple state transfer}, where there is a subset of vertices such that every pair of vertices in that subset has perfect state transfer between them, and are particularly interested in the non-trivial cases where multiple state transfer occurs between three or more vertices. The oriented \( K_3 \) is still an example of graphs with multiple state transfer, but it is not the only one.

\begin{figure}
  \begin{center}
    \begin{tikzpicture}[node distance = 2 cm]
      \tikzset{VertexStyle/.style = {shape= circle, ball color= blue!20!white, inner sep= 2pt, outer sep= 0pt, minimum size = 24 pt}}
      \tikzset{EdgeStyle/.style   = {->, double= blue, double distance = 1pt}}
    
      \draw (315:1 cm) node[VertexStyle, font=\small] (1) {1} ;
      \draw (135:1 cm) node[VertexStyle, font=\small] (6) {6} ;
      \draw (225:1 cm) node[VertexStyle, font=\small] (2) {2} ;
      \draw (45:3 cm) node[VertexStyle, font=\small] (3) {3} ;

      \draw (45:1 cm) node[VertexStyle, font=\small] (5) {5} ;
      \draw (135:3 cm) node[VertexStyle, font=\small] (0) {0} ;
      \draw (315:3 cm) node[VertexStyle, font=\small] (7) {7} ;
      \draw (225:3 cm) node[VertexStyle, font=\small] (4) {4} ;
      
      \tikzset{VertexStyle/.style = {shape= circle, ball color= white, inner sep= 0pt, outer sep= 0pt, minimum size = 0 pt}}
    
      \draw (315:4 cm) node[VertexStyle, font=\small] (8) {} ;
      \draw (225:4 cm) node[VertexStyle, font=\small] (9) {} ;
      
      \draw [EdgeStyle] (0) edge (2) ;
      \draw [EdgeStyle] (0) edge (3) ;
      \draw [EdgeStyle] (0) edge (4) ;
      \draw [EdgeStyle] (0) edge (5) ;
      \draw [EdgeStyle] (6) edge (0) ;
      \draw [EdgeStyle] (9) edge [bend right=30] (7) ;

      \draw [EdgeStyle] (1) edge (2) ;
      \draw [EdgeStyle] (1) edge (3) ;
      \draw [EdgeStyle] (1) edge (4) ;
      \draw [EdgeStyle] (1) edge (5) ;
      \draw [EdgeStyle] (1) edge (6) ;
      \draw [EdgeStyle] (7) edge (1) ;

      \draw [EdgeStyle] (6) edge (2) ;
      \draw [EdgeStyle] (6) edge (3) ;
      \draw [EdgeStyle] (6) edge (4) ;
      \draw [EdgeStyle] (6) edge (5) ;
    
      \draw [EdgeStyle] (7) edge (2) ;
      \draw [EdgeStyle] (7) edge (3) ;
      \draw [EdgeStyle] (7) edge (4) ;
      \draw [EdgeStyle] (7) edge (5) ;
      
      \draw [EdgeStyle] (4) edge (2) ;
      \draw [EdgeStyle] (2) edge (5) ;

      \draw [EdgeStyle] (8) edge [ bend left=30] (4) ;
      \draw [EdgeStyle] (5) edge (3) ;
      
      \tikzset{EdgeStyle/.style   = {double= blue, double distance = 1pt}}
      \draw [EdgeStyle] (0) edge [bend right=60] (9) ;
      \draw [EdgeStyle] (3) edge [bend left=60] (8) ;
    \end{tikzpicture}
  \end{center}
  \caption{A New Example of Multiple State Transfer} \label{multi8} 
\end{figure}
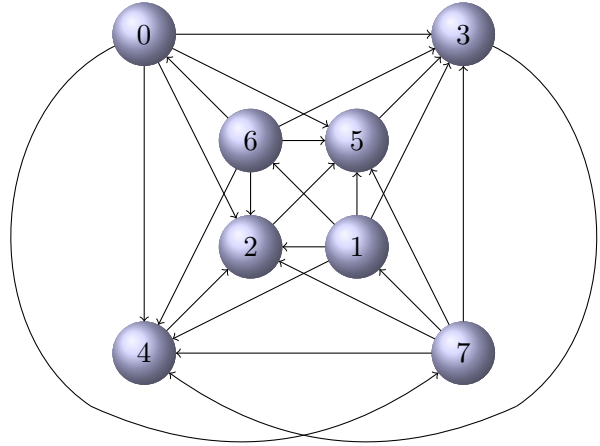

The graph in Figure~\ref{multi8} has multiple state transfer on the set of vertices \( \lsb 0, 1, 6, 7 \rsb. \) In both this example and in the oriented \( K_3, \) there is a switching automorphism permuting the vertices involved in multiple state transfer. In this case, determining when multiple state transfer occurs simplifies to determining when perfect state transfer occurs between two vertices determined by the automorphism.

\begin{lem}\label{mst}
	Let \( X \) be an oriented graph with switching automorphism \( \tilde{P} \) of order \( n \). If perfect state transfer occurs 
	from vertex \( a \) to vertex \( \tilde{P} a \), then there is multiple state transfer on
 	\[ \lsb a, \tilde{P}a, \ldots, \tilde{P}^{n-1} a \rsb. \]
\end{lem}

\proof 
Let \( \tilde{P}^k a, \tilde{P}^{\ell} a\) be two vertices in the set. By Corollary~\ref{multipleAuto}, perfect state transfer occurs from \( a \) to \( \tilde{P}^k a \) at time \( k \tau \) and from \( a \) to \( \tilde{P}^{\ell} a \) at time \( \ell \tau \). Then
\begin{align*}
  U((\ell - k) \tau) \tilde{P}^k \ver{a} &= U(\ell \tau) U( -k \tau) \tilde{P}^k \ver{a} \\
  &= \pm U(\ell \tau)\ver{a} \\
  &= \pm \tilde{P}^{\ell} \ver{a}.
\end{align*}

Thus there is perfect state transfer from \( \tilde{P}^k a \) to \( \tilde{P}^{\ell} a \) and, 
since \( k \) and \( \ell \) were arbitrary, there is multiple state transfer on the set.\qed

For a switching automorphism \( \tilde{P} \) and a vertex \( a \), we may adapt the tools used to study perfect state transfer on unoriented graphs to develop a characterization for when perfect state transfer occurs from \( a \) to \( \tilde{P} a \). This allows us to come up with a characterization of multiple state transfer on the image under a switching automorphism.

\section{Cospectrality}

Let \( a \) be a vertex. We define the \textit{eigenvalue support} at \( a, \) denoted \( \Phi_a \), to be the set
\[ \lsb \theta_r : E_r \ver{a} \neq \mathbf{0} \rsb. \]

Two vertices \( a \) and \( b \) are \textit{strongly cospectral} if for each spectral idempotent \( E_r \) 
with \( \theta_r \in \Phi_a \), there exists some complex number \( \alpha_r \) with \( \abs{\alpha_r} = 1 \) 
such that
\[ E_r \ver{a} = \alpha_r E_r \ver{b}. \]
More details about strongly cospectral vertices can be found in Godsil and Smith~\cite{Godsil2017e}. For unoriented graphs, the entries of the spectral idempotents will be real and thus \( \alpha = \pm 1 \). Additionally, if there is perfect state transfer between two vertices, they are strongly cospectral and strong cospectrality is one of several necessary and sufficient conditions for perfect state transfer to occur, as shown by Coutinho in Chapter~2 of his thesis~\cite{Coutinho2014}. 

\begin{lem}
	\label{lem:pst-scsp}
	Let \( X \) be an oriented graph with vertices \( a \) and \( b \). If there is perfect state transfer from \( a \) to \( b \), then \( a \) and \( b \) 
	are strongly cospectral.
\end{lem}

\proof 
Suppose there is perfect state transfer from \( a \) to \( b \) at time \( \tau. \) Then for all eigenvalues \( \theta_r \), we may define the corresponding \( \alpha_r \) to be either \( e^{-\tau \theta_r} \) or \( - e^{-\tau \theta_r} \) depending on whether the phase factor is plus or minus one.

Because there is perfect state transfer from \( a \) to \( b \), we know that
\[ \sum_r e^{\tau \theta_r} E_r \ver{a} = \pm \sum_r E_r \ver{b}, \]
or equivalently, for all \( r \),
\[ e^{\tau \theta_r} E_r \ver{a} = \pm E_r \ver{b}. \]
By multiplying by \( e^{- \tau \theta_r} \) we get
\[ E_r \ver{a} = \pm e^{- \tau \theta_r} E_r \ver{b} = \alpha_r E_r \ver{b}, \]
and therefore \( a \) and \( b \) are strongly cospectral.\qed


If we have perfect state transfer from $a$ to $b$, it follows from the discussion above
that there is a real number $q_r(a,b)$ such that 
\[
	-1 \le q_r(a,b) \le 1
\]
and
\[
	E_r\ver{a} = e^{i\pi q_r(a,b)} E_r\ver{b}.
\]
We refer to $q_r(a,b)$ as the quarrel from $a$ to $b$ (relative to the eigenvalue $\theta_r$.)

\begin{lem}\label{multipleRobust} 
  Let \( X \) be an oriented graph with switching automorphism \( \tilde{P}. \) If vertices \( a \) and \( \tilde{P} a \) are 
	strongly cospectral then, for any positive integer \( n \), vertices \( a \) and \( \tilde{P}^n a \) 
	are strongly cospectral, and
  	\[ e^{i \pi q_r(a, \tilde{P}^n a)} = e^{i \pi n q_r(a, \tilde{P} a)}. \]
\end{lem}

\proof 
We will prove this by induction. When \( n = 1, \) we have
\[ E_r \ver{a} = e^{i \pi q_r(a,b)} E_r \tilde{P} \ver{a}, \]
which is the definition. Let \( k \) be a positive integer such that \( a \) and \( \tilde{P}^k a \) are strongly cospectral, and let \( q = q_r \lp a, \tilde{P} a \rp. \) Then
\begin{align*}
  E_r \ver{a} &= e^{i \pi k q} E_r \tilde{P}^k \ver{a} \\
  &= e^{i \pi k q} \tilde{P}^k \lp \tilde{P}^T \rp^k E_r \tilde{P}^k \ver{a} \\
  &= e^{i \pi k q} \tilde{P}^k E_r \ver{a} \\
  &= e^{i \pi \lp k+1 \rp q} \tilde{P}^k E_r \tilde{P} \ver{a} \\
  &= e^{i \pi \lp k+1 \rp q} \tilde{P}^k E_r \lp \tilde{P}^T \rp^k \tilde{P}^{k+1} \ver{a} \\
  &= e^{i \pi \lp k+1 \rp q} E_r \tilde{P}^{k+1} \ver{a}.
\end{align*}
Note that
\[ \ver{a}^T E_r \ver{a} = e^{i \pi \lp k+1 \rp q} \ver{a}^T E_r \tilde{P}^{k+1} \ver{a} \]
is a positive real number, and thus
\[ e^{i \pi q_r \lp a, \tilde{P}^{k+1} a \rp} = e^{i \pi \lp k+1 \rp q} \]
and therefore \( a \) and \( \tilde{P}^{k+1} a\) are strongly cospectral with the desired condition on the quarrel.\qed

We now offer a first characterization of when there will be perfect state transfer between two vertices. This characterization applies for all graphs, and is an oriented analogue of Theorem 2.4.2 in~\cite{Coutinho2014}.

\begin{thm}\label{firstChar}Let \( X \) be an oriented graph with vertices \( a \) and \( b \). There is perfect state transfer from \( a \) to \( b \) if and only if:
  \begin{enumerate}[label = (\roman*)]
  \item Vertices \( a \) and \( b \) are strongly cospectral.
  \item There exists some \( \tau \in \mathbb{R} \) such that, for all \( \theta_r \in \Phi_a \),
    \[ q_r(a,b) + \frac{i \tau \theta_r}{\pi} \]
    is always an even integer or always an odd integer.
  \end{enumerate}
\end{thm}

\proof By Lemma~\ref{lem:pst-scsp}, we may assume \( a \) and \( b \) are strongly cospectral and show that perfect state transfer from \( a \) to \( b \) occurs if and only if (ii) holds. Perfect state transfer occurs if and only if there exists some \( \tau \in \mathbb{R} \) and some phase factor \( \pm 1 \) such that
\[ e^{\tau \theta_r} E_r \ver{a} = \pm E_r \ver{b} = e^{i \pi q_r(a,b)} E_r \ver{a}. \]

If the phase factor is \( 1\), this is equivalent to saying that for all \( \theta_r \in \Phi_a \),
\[ e^{ \tau \theta_r - i \pi q_r(a,b)} = 1, \]
or there exists some integer \( k_r \) such that
\[ \tau \theta_r - i \pi q_r(a,b) = 2 k_r i \pi. \]
Multiplying through by \( \frac{i}{\pi} \) shows us that perfect state transfer occurs at time \( \tau \) with phase factor \( 1 \) if and only if
\[ \frac{i \tau \theta_r}{\pi} + q_r(a,b) \]
is always an even integer.

Similarly, if the phase factor is \( -1 \), then perfect state from \( a \) to \( b \) at time \( \tau \) is equivalent to saying that for all \( \theta_r \in \Phi_a \), there exists some integer \( k_r \) such that
\[ \tau \theta_r - i \pi q_r(a,b) = \lp 2 k_r +1 \rp i \pi, \]
or equivalently
\[ \frac{i \tau \theta_r}{\pi} + q_r(a,b) \]
is always an odd integer.\qed

This provides us with a first characterization, which applies for all oriented graphs. However, it 
relies upon knowing at what time perfect state transfer occurs. In order to deal with this, we consider a 
closely related but easier to characterize property of quantum walks.

\section{Periodicity}

A vertex \( a \) in an oriented graph is \textit{periodic} if there is a time \(\sigma\) such that
\(U(\sigma)\ver{a}= \pm \ver{a}\). (For general graphs, a vertex is periodic if there is time \(\sigma\)
and a complex scalar \(\gamma\) with norm \(1\) such that \(U(\sigma)\ver{a}=\gamma\ver{a}\),
further if there is perfect state transfer from \(a\) to \(b\), then both vertices \(a\)
and \(b\) are periodic. For more detail, see \cite{Godsil2017}.)

For oriented graphs with a switching automorphism of order \( n \) taking \( a \) to \( b \) and perfect 
state transfer from \( a \) to \( b \) at time \( \tau \), Corollary~\ref{multipleAuto} gives a
proof that the vertex \( a \) will be periodic at time \( n \tau. \)

Periodicity is a simpler property than perfect state transfer, and it is much easier to characterize
for both oriented and unoriented graphs. For unoriented graphs, Godsil~\cite{GodsilWhen} characterized graphs with periodic 
vertices based only on the elements of the eigenvalue support. Here, we prove 
a similar, but simpler characterization for oriented graphs.

\begin{thm}\label{periodic}Let \( X \) be a connected oriented graph with at least two vertices. Then the following are equivalent:
\begin{enumerate}[label = (\roman*)]
\item The vertex \( a \) is periodic.
\item For all \( r, s \) with \( \theta_r, \theta_s \in \Phi_a \) and \( \theta_s \neq 0 \), the ratio \( \frac{\theta_r}{\theta_s} \) is rational.
\item There exists a square-free positive integer \( \Delta \) such that all eigenvalues in \( \Phi_a \) are in \( \mathbb{Z}(\sqrt{ -\Delta}) \).
\end{enumerate}
\end{thm}

\proof Suppose that \( a \) is a periodic vertex. This is equivalent to saying that for all \( \theta_{r} \in \Phi_a, \)
\[ e^{\tau \theta_{r}} E_{r} \ver{a} = \pm E_{r} \ver{a}. \]

In particular, for any \( \theta_r, \theta_s \in \Phi_a \) where \( \theta_s \neq 0, \) this means that \( e^{\tau \theta_r} = e^{\tau \theta_s} = \pm 1 \), and therefore \( \tau \theta_r \) and \( \tau \theta_s \) are both integer multiples of \( i \pi \). It follows that
\[ \frac{\theta_r}{\theta_s} = \frac{\tau \theta_r}{\tau \theta_s} \in \mathbb{Q}. \]

Next, suppose that for all \( r, s \) with \( \theta_r \in \Phi_a \) and \( \theta_s \neq 0 \), the ratio \( \frac{\theta_r}{\theta_s} \) is rational. Since \( X \) is connected, we know that
\[ \sum_r \theta_r E_r \ver{a} = A \ver{a} \neq \mathbf{0}, \]
and in particular there must be some \( \theta_0 \neq 0 \) in \( \Phi_a \). For all \( \theta_r \in \Phi_a \), we know there exists some rational number \( m_r \) such that \( \theta_r = m_r \theta_0 \) and thus
\[ \prod_{\theta_r \in \Phi_a} \theta_r = \theta_0^{\abs{\Phi_a}} \prod_{\theta_r \in \Phi_a} m_r. \]

The eigenvalue support is closed under algebraic conjugates, so any automorphism of the splitting field of the characteristic polynomial of the adjacency matrix fixes
\[ \prod_{\theta_r \in \Phi_a} \theta_r. \]

It follows from field theory that \( \theta_0^{\abs{\Phi_a}} \in \mathbb{Q} \), and we may let \( d \) be the minimal integer such that \( \theta_0^d \in \mathbb{Q} \). Then \( \theta_0 \) has conjugate eigenvalues that differ by \( d \)th roots of unity and, since the eigenvalues of skew symmetric matrices have no real part, \( d \leq 2 \) and \( \theta_0^2 \in \mathbb{Q}. \) Because \( \theta_0 \) is an eigenvalue of a matrix with integer entries and therefore an algebraic integer, and rational algebraic integers are integers, we have in fact that \( \theta_0^2 \in \mathbb{Z} \). Letting \( \Delta \) be the largest square-free integer dividing \( \theta_0^2 \), we write \( \theta_0 = a \sqrt{-\Delta} \) for some integer \( a \).

Then for any \( \theta_r \in \Phi_a \), we have
\[ \theta_r = a m_r \sqrt{-\Delta} \]
and
\[ \theta_r^2 = - (a m_r)^2 \Delta \in \mathbb{Z}. \]
Since \( \Delta \) is square-free, it must be the case that \( a m_r \in \mathbb{Z} \), and therefore all eigenvalues in \( \Phi_a \) are in \( \mathbb{Z} (\sqrt{-\Delta}). \)

Finally, suppose that there exists a square-free positive integer \( \Delta \) such that all eigenvalues in \( \Phi_a \) are in \( \mathbb{Z} (\sqrt{ -\Delta}) \). Then
\[ U(\frac{2 \pi}{\sqrt{\Delta}}) \ver{a} = \sum_r e^{\frac{2 \pi \theta_r}{\sqrt{\Delta}}} E_r \ver{a} = \sum_{\theta_r \in \Phi_a} E_r \ver{a} = \ver{a}, \]
and so \( a \) is periodic.\qed

From this characterization, we see that the eigenvalues are sufficient to determine when a vertex is periodic. They can also tell us the \textit{minimum period}, or the first time that a vertex is periodic.

\begin{lem}\label{minPeriod}
  Let \( X \) be a connected oriented graph with periodic vertex \( a, \) and let \( \Delta \) be the square-free positive integer such that for all \( \theta_r \in \Phi_a, \) \( \theta_r \in \mathbb{Z} \lb \sqrt{-\Delta} \rb. \) Assume
  \[ g = \mathrm{ gcd} \left(\lsb \frac{\theta_r}{\sqrt{-\Delta}} \rsb_{\theta_r \in \Phi_a}\right). \]
  Then if the phase factor of \(a\) is \( -1 \), the minimum period is
  \[ \frac{\pi}{g \sqrt{\Delta}} \]
  and if the phase factor is \( 1 \), the minimum period is
  \[ \frac{2 \pi}{g \sqrt{\Delta}} \]
\end{lem}

\proof 
Suppose the minimum period is \( \sigma \) with phase factor \( -1 \). Then for all \( \theta_r \in \Phi_a \)
\[ e^{\sigma \theta_r} = -1 \]
so
\[ \sigma \theta_r = i \pi k_r \]
for some odd integer \( k_r \). But \( \theta_r \) is an imaginary quadratic integer, and so the smallest that \( \sigma \) can be is \( \frac{\pi}{g \sqrt{\Delta}} \).

Similarly, if the minimum period \( \sigma \) has phase factor \( 1 \) then for all \( \theta_r \in \Phi_a \)
\[ e^{\sigma \theta_r} = 1 \]
so
\[ \sigma \theta_r = 2 i \pi k_r \]
for some integer \( k_r \), so the smallest that \( \sigma \) can be is \( \frac{2 \pi}{g \sqrt{\Delta}} \).\qed

\section{Characterization}

Using our earlier characterization and what we know about periodicity, we can come up with a better characterization of when a vertex has perfect state transfer to its image under a switching automorphism, similar to Theorem 2.4.4 for nonoriented graphs~\cite{Coutinho2014}. This gives a characterization of when multiple state transfer can occur that relies only the spectral decomposition of the graph, a given switching automorphism, and a vertex in the graph. 

\begin{thm}\label{completeChar}
	Let \( X \) be an oriented graph with vertex \( a \), 
	and let \( \tilde{P} \) be a switching automorphism of order \( n \). Then there is multiple 
	state transfer on the set \( \{a, \tilde{P} a, \ldots, \tilde{P}^{n-1} a\} \) if and only if:
  \begin{enumerate}[label = (\roman*)]
  \item Vertices \( a \) and \( \tilde{P}a \) are strongly cospectral.
  \item Vertex \( a \) is periodic.
  \item Let
    \[ g = \mathrm{ gcd} \lp \lsb \frac{\theta_r}{\sqrt{-\Delta}} \rsb_{\theta_r \in \Phi_a} \rp. \]
    Then either
    \begin{enumerate}[label = (\alph*)]
    \item Vertex \( a \) has periodic phase factor \( -1 \) and there exists an integer \( m < n \) 
	relatively prime to \( n \) such that
      \[ \frac{\theta_r}{n g \sqrt{-\Delta}} + m q_r(a, \tilde{P} a) \]
      is always an odd integer or always an even integer for all \( \theta_r \in \Phi_a. \)
    \item Vertex \( a \) has periodic phase factor \( 1 \) and there exists an integer \( m < n \) 
	relatively prime to \( n \) such that
      \[ \frac{2 \theta_r}{n g \sqrt{-\Delta}} + m q_r(a, \tilde{P} a) \]
      is always an odd integer or always an even integer for all \( \theta_r \in \Phi_a. \)
    \end{enumerate}
  \end{enumerate}
\end{thm}

\proof 
If there is multiple state transfer, there is necessarily perfect state transfer from \( a \) 
to \( \tilde{P} a\). Then by Lemma~\ref{lem:pst-scsp} and Lemma~\ref{periodic} the first two conditions hold, 
so we may assume them and prove that multiple state transfer occurs on the set if and only if (iii) holds. 
Since \( a \) is periodic, it can either have phase factor \( - 1 \) or 1. 
Letting \( \sigma \) denote the minimum period, by Lemma~\ref{minPeriod}, (iii) is equivalent to saying that 
for some \( m \)
\[ \frac{\sigma \theta_r}{\pi n i} + m q_r(a, \tilde{P} a) \]
is always an even integer or always an odd integer.

Suppose (iii) holds. By Lemma~\ref{multipleRobust}, vertices \( a \) and \( \tilde{P}^m a \) are strongly cospectral and there exists some integer \( k \) such that
\[ q_r(a, \tilde{P}^m a) = m q_r(a, \tilde{P} a) - 2k. \] 
Then
\[ \frac{\sigma \theta_r}{\pi n i} + m q_r(a, \tilde{P} a) 
	= \frac{\sigma \theta_r}{\pi n i} -2k + q_r(a, \tilde{P}^m a) \]
is always an even integer or always an odd integer, so by Theorem~\ref{firstChar} there is perfect 
state transfer from \( a \) and \( \tilde{P}^m a \) at \( \frac{\sigma}{n \pi} \). Then by 
Lemma~\ref{mst}, since \( m \) is relatively prime to \( n \) there is multiple state transfer on the set
\[ \{ a, \tilde{P}^m a, \ldots, \tilde{P}^{m (n-1)} a \} 
	= \{a, \tilde{P} a, \ldots, \tilde{P}^{n-1} a. \} \]

Conversely, suppose that multiple state transfer occurs on this set and for \( 1 \leq i \leq n-1 \) let \( \tau_i \) be the first time that perfect state transfer occurs from \( a \) to \( \tilde{P}^i a \). Note that these will all be distinct times strictly less than \( \sigma \) and that therefore \( n \tau_i \) is a distinct period for each \( i \). Every period is an integer multiples of the minimum period, therefore \( n \tau_i \) is some integer multiple of \( \sigma \), and thus the set of minimal perfect state transfer times from \( a \) are
\[ \lsb \frac{z_1 \sigma}{n}, \ldots, \frac{z_{n-1} \sigma}{n} \rsb \]
where \( z_1, \ldots, z_{n-1} \) are \( n-1 \) distinct integers between 1 and \( n - 1 \). Therefore there must be some \( m \) such that \( \tau_m = \frac{\sigma}{n}. \) Since \( \sigma \) is the minimum period, \( m \) must be relatively prime to \( n \) or else \( \tilde{P}^m \) would be a nontrivial switching automorphism of order \( m' < m \) and thereforem \( m' \tau_m \) would be a period strictly less than \( \sigma. \)

By Lemma~\ref{multipleRobust}, vertices \( a \) and \( \tilde{P}^m a\) are strongly cospectral and there exists some integer \( k \) such that
\[ q_r(a, \tilde{P}^m a) = 2k + m q_r(a, \tilde{P} a). \]
Therefore by Theorem~\ref{firstChar}, for all \( \theta_r \in \Phi_a \)
\[ q_r(a, \tilde{P}^m) + \frac{i \tau_m \theta_r}{\pi} 
	= m q_r(a, b) + 2k + \frac{i \sigma \theta_r}{n \pi} \]
is always an even integer or always an odd integer, and thus multiple state transfer implies (iii).\qed

\section{Further Questions}

\begin{figure}
  \begin{center}
    \begin{tikzpicture}[node distance = 2 cm]
      \tikzset{VertexStyle/.style = {shape= circle, ball color= blue!20!white, inner sep= 2pt, outer sep= 0pt, minimum size = 24 pt}}
      \tikzset{EdgeStyle/.style   = {->, double= blue, double distance = 1pt}}
      
      \draw (90:2.25 cm) node[VertexStyle, font=\small] (0) {0} ;
      \draw (210:2.75 cm) node[VertexStyle, font=\small] (1) {1} ;
      \draw (330:2.75 cm) node[VertexStyle, font=\small] (2) {2} ;
    
      \draw (180:.65 cm) node[VertexStyle, font=\small] (3) {3} ;
      \draw (0:.65 cm) node[VertexStyle, font=\small] (4) {4} ;
    
      \tikzset{EdgeStyle/.style   = {->, double= blue, double distance = 1pt}}
      \draw [EdgeStyle] (3) edge (4) ;
    
      \draw [EdgeStyle] (0) edge [bend right] (3) ;
      \draw [EdgeStyle] (0) edge [bend left] (4) ;
      \draw [EdgeStyle] (1) edge [bend left] (3) ;
      \draw [EdgeStyle] (1) edge [bend right] (4) ;
      \draw [EdgeStyle] (2) edge [bend left] (3) ;
      \draw [EdgeStyle] (2) edge [bend right] (4) ;
    
    \end{tikzpicture}
  \end{center}
 \caption{Perfect state transfer occurs at an irrational multiple of the period}  \label{ununderstandable} 
\end{figure}
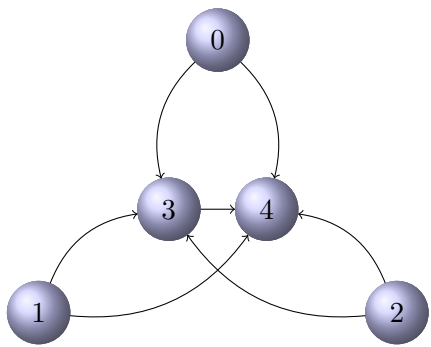

The characterization in Theorem~\ref{completeChar} applies whenever there is perfect state transfer and a switching automorphism between two vertices. In this case, the time for perfect state transfer will be a rational multiple of the minimum period, which is itself a rational multiple of \( \frac{\pi}{\sqrt{\Delta}} \). However, this is not the only situation under which perfect state transfer can occur. The graph in Figure~\ref{ununderstandable} has no switching automorphism between vertices 3 and 4, since the other three vertices are oriented towards vertex 3  in the same way that they are oriented towards vertex 4, but there is a single edge oriented from vertex 3 to vertex 4. However, Theorem~\ref{firstChar} can be used to verify that perfect state transfer occurs from vertex 3 to vertex 4 at time\[ \frac{\pi - \arccos\lp\frac{3}{4} \rp}{\sqrt{7}}, \]
but \( \arccos \lp \frac{3}{4} \rp \) is not a rational multiple of \( \pi \)~\cite{Lato2019}. This example raises several questions:
\begin{itemize}
\item Is it possible to improve the characterization in Theorem~\ref{firstChar} so that it does not depend on the time \( \tau \) while still having it apply to all graphs?
\item What is the relationship between the period and the first time to perfect state transfer if a graph does not have a switching automorphism between the vertices involved in perfect state transfer?
\item Are there examples of multiple state transfer that do not arise from switching automorphisms?
\end{itemize}

\end{document}